\documentclass[12pt]{amsart}
\usepackage{amscd,amssymb}
\usepackage[graph,frame,poly,arc]{xy}  
\usepackage[plainpages,backref,urlcolor=blue]{hyperref}
\usepackage{mathtools}

\theoremstyle{plain}
\newtheorem{thm}[subsection]{Theorem}

\newtheorem{prop}[subsection]{Proposition}
\newtheorem{cor}[subsection]{Corollary}

\theoremstyle{definition}
\newtheorem{rk}[subsection]{Remark}
\newtheorem{definition}[subsection]{Definition}

\newtheorem{conj}[subsection]{Conjecture}
\newtheorem{question}[subsection]{Question}

\numberwithin{equation}{section}
\newcommand{\OO}{{\mathcal O}}

\newcommand{\A}{{\mathcal A}}
\newcommand{\B}{{\mathcal B}}

\newcommand{\C}{\mathbb{C}}
\newcommand{\PP}{\mathbb{P}}

\newcommand{\dd}{{\rm d}}

\DeclareMathOperator{\im}{im}

\DeclareMathOperator{\mult}{mult}
\DeclareMathOperator{\defect}{def}

\DeclarePairedDelimiter\ceil{\lceil}{\rceil}

\begin{document}
%\date{June 4, 2009}

\title [On free curves and related open problems]
{On free curves and related open problems}

\author[Alexandru Dimca]{Alexandru Dimca}
\address{Universit\'e C\^ ote d'Azur, CNRS, LJAD, France and Simion Stoilow Institute of Mathematics,
P.O. Box 1-764, RO-014700 Bucharest, Romania}
\email{dimca@unice.fr}



\begin{abstract} 
In this paper we collect the main properties of free curves in the complex projective plane and a lot of conjectures and open problems, both old and new. In the quest to understand the mystery of free curves, many tools were developed and many results were obtained, which apply to any reduced plane curve,
and some of them are recorded here.
\end{abstract}
 
\maketitle

%\tableofcontents

\section{Introduction} 

The notion of free hypersurface, or free divisor, has a long and fascinating history. This concept was considered mostly in the following two settings, namely
\begin{enumerate}

\item[{(L)}] The local setting, where one looks at germs of hypersurfaces, as in the seminal paper by K. Saito \cite{KS}. In particular, any plane curve singularity is known to be free.

\item[{(P)}]  The projective setting, where one considers hypersurfaces in a complex projective space $\PP^n$, as for instance in \cite{B+,Dmax,MRL,Sim2,ST,To}. In particular, the free projective hypersurfaces satisfy a lot of restrictions, for instance they must be rather singular, namely the singular locus must have codimension one.
\end{enumerate}

Note that the  central hyperplane arrangements in some affine space $\C^m$ can be considered from both view-points, and Terao's Conjecture, see Section 3 below for the case $m=3$, is a beautiful open problem in this area.
More generally, any hypersurface in $\C^m$ defined by a homogeneous polynomial can be considered from both view-points, and freeness in one setting is the same as freeness in the other.

In this paper we consider only reduced curves in the complex projective plane $\PP^2$. In Section 2 we recall some of the main properties of free curves
as well as a number of general results on plane curves, helping us in the understanding of free curves.

In Section 3 we discuss several stronger versions of Terao's Conjecture
in the case of line arrangements in $\PP^2$. Though a lot of effort was devoted to this conjecture over the last 40 years or so, only partial results are known, see \cite{Abe, BK, DHA,OT,To,Yo14}. The main interest of these stronger versions of Terao's Conjecture is that {\it they refer to properties of any line arrangement, and hence they tell something new about any line arrangement, and sometimes about any plane curve}. In Remark \ref{rkZ} we explain how starting with a pair of line arrangements $\A$ and $\A'$, having the same intersection lattice but distinct invariants related to the minimal resolution of their Jacobian ideal, for instance the arrangements constructed by Ziegler in \cite{Z},
new such pairs can be obtained by adding new lines.

In Section 4 we discuss new invariants, related to the defects of some linear systems with respect to the singular subscheme of a plane curve, and state a new equivalent version of Terao's Conjecture, see Conjecture \ref{conj00}, and dismiss a natural
stronger version of Terao's Conjecture, see Question \ref{qT4}.

In Section 5 we discuss the relation between rational cuspidal plane curves and the free curves. It is the study of this relation that conducted us to introduce the notion of nearly free curve in \cite{DStRIMS}. The conjecture, saying that any rational cuspidal plane curve is either free or nearly free is known to hold for all even degree curves and in 'most' cases for odd degree curves. However, this conjecture is still open in full generality.

In Section 6 we start by recalling a result due to A. du Plessis and C.T.C. Wall, see Theorem \ref{thmCTC}, which gives lower and upper bounds for the total Tjurina number $\tau(C)$ of a reduced curve $C$. In some cases, the curves realizing the maximal value for $\tau(C)$ are exactly the free curves. In this section we investigate the non-free curves for which these bounds are attained and put forth a number of  open questions.

In Section 7 we discuss the existence of some free curves having only simple singularities of type $A,D,E$. These curves are strongly related to the maximizing curves introduced by U. Persson \cite{Persson}. For such a maximizing curve $C$, the associated smooth surface $\widetilde{X}$ obtained as the minimal resolution of the double cover $X$ of $\PP^2$ ramified along $C$ might have the maximal possible Picard number. For this geometric construction, the degree of $C$ must be even. However, the odd degree curves with similar properties from the view-point of free curves were introduced in \cite{DPok} and until now their existence is limited to odd degrees up to 9.

Finally in Section 8 we introduce the supersolvable plane curves in analogy to supersolvable line arrangements, a classical notion in hyperplane arrangement theory. We conjecture that such supersolvable curves are free in Conjecture \ref{conj18}. Then we recall both a general result supporting this conjecture, see Theorem \ref{thm58}, and several, in fact countable many examples of supersolvable curves, where the freeness is not a consequence of this result.

The author would like to thank all the collaborators with which he worked over the years on this beautiful and highly mysterious subject. Most of the results quoted in this very personal survey are the result of our joint work.

\section{Basic facts and notions related to free curves } 

Let $S=\C[x,y,z]$ be the polynomial ring in three variables $x,y,z$ with complex coefficients, and let $C:f=0$ be a reduced curve of degree $d\geq 3$ in the complex projective plane $\PP^2$. 
We denote by $J_f$ the Jacobian ideal of $f$, i.e. the homogeneous ideal in $S$ spanned by the partial derivatives $f_x,f_y,f_z$ of $f$, and  by $M(f)=S/J_f$ the corresponding graded quotient ring, called the Jacobian (or Milnor) algebra of $f$.
Consider the graded $S$-module of Jacobian syzygies of $f$ or, equivalently, the module of derivations killing $f$, namely
\begin{equation}
\label{eqD0}
D_0(f)=\{\rho=(a,b,c) \in S^3 \ : \ af_x+bf_y+cf_z=0\}.
\end{equation}
This module is also denoted in the literature by 
$AR(f)$ (all Jacobian relations for $f$) and ${\rm Syz}(f)$, the Jacobian syzygies of $f$.
According to Hilbert Syzygy Theorem, the graded Jacobian algebra $M(f)$ has a minimal free resolution of the form
\begin{equation}
\label{res1}
0 \to F_3 \to F_2 \to F_1 \to F_0,
\end{equation}
where clearly $F_0=S$, $F_1=S^3(1-d)$ and the morphism $F_1 \to F_0$ is given by
$$(a,b,c) \mapsto af_x+bf_y+cf_z.$$
With this notation, the graded $S$-module of Jacobian syzygies $D_0(f)$ has the following minimal resolution
$$0 \to F_3(d-1) \to F_2(d-1).$$
We say that $C:f=0$ is an {\it $m$-syzygy curve} if the module  $F_2$ has rank $m$. Then the module $D_0(f)$ is generated by $m$ homogeneous syzygies, say $\rho_1,\rho_2,...,\rho_m$, of degrees $d_j=\deg r_j$ ordered such that $$1 \leq d_1\leq d_2 \leq ...\leq d_m.$$ 
We call these degrees $(d_1, \ldots, d_m)$ the {\it exponents} of the curve $C$ and $\rho_1,...,\rho_m$ a {\it minimal set of generators } for the module  $D_0(f)$. 
The smallest degree $d_1$ is sometimes denoted by ${\rm mdr}(f)$ and is called the minimal degree of a Jacobian relation for $f$. 

The curve $C$ is {\it free} when $m=2$, since then  $D_0(f)$ is a free module of rank 2, see for instance \cite{KS,Sim2,ST,To}. Moreover, there are two classes of 3-syzygy curves which are intensely studied, since they are in some sense the closest to free curves.
First, we have the  {\it nearly free curves}, introduced in \cite{DStRIMS} and studied in \cite{AD, B+, Dmax,  MaVa} which are 3-syzygy curves satisfying $d_3=d_2$ and $d_1+d_2=d$. 
Then, we have the  {\it plus-one generated line arrangements} of level $d_3$, introduced by Takuro Abe in \cite{Abe} and recently studied in \cite{MP,MV},
which are 3-syzygy line arrangements satisfying $d_1+d_2=d$. In general, a  3-syzygy curve will be called a {\it plus-one generated curve}  if it satisfies $d_1+d_2=d$. { In particular, nearly free curves are a special type of plus-one generated curves.} The study of these classes of curves goes naturally together, as the following result shows, see for it and other similar results \cite{POG,MP2}.

\begin{thm}
\label{thm1}
Let $C:f=0$ be a reduced curve in $\PP^2$, $L$ a line in $\PP^2$, which is not an irreducible component of $C$. We assume that the union $C'=C \cup L:f'=0$ is a free curve. Then the curve $C$ is either free or a plus-one generated curve.
\end{thm}

We have the following characterizations for free and plus-one generated curves, see  \cite[Theorem 2.3]{minTjurina}, where $D_0(f)$ is denoted by $AR(f)$.

\begin{thm}
\label{thmA}
Let $C:f=0$ be a reduced plane curve of degree $d$ and let $d_1$ and $d_2$ be the minimal degrees of a minimal system of generators for the module of Jacobian syzygies $D_0(f)$ as above.
Then the following hold.
\begin{enumerate}

\item The curve $C$ is free if and only if $d_1+d_2=d-1$.

\item The curve $C$ is plus-one generated if and only if $d_1+d_2=d$.

\item In all the other cases $d_1+d_2 >d$.

\end{enumerate}

\end{thm}

Let $I_f$ denote the saturation of the ideal $J_f$ with respect to the maximal ideal ${\bf m}=(x,y,z)$ in $S$ and consider the following  local cohomology group, usually called the Jacobian module of $f$, 
 $$N(f)=I_f/J_f=H^0_{\bf m}(M(f)).$$
We set $n(f)_k=\dim N(f)_k$ for any integer $k$ and introduce the {\it freeness defect of the curve} $C$ by the formula
$$\nu(C)=\max _j \{n(f)_j\}$$ as in \cite{AD}.
Note that $C$ is free if and only if $N(f)=0$, see for instance \cite{ST}, and hence in this case $\nu(C)=0$, and $C$ is nearly free if and only if $\nu(C)=1$, see \cite{DStRIMS}. If we set $T=3(d-2)$, then the sequence $n(f)_k$ is symmetric with respect to the middle point $T/2$, that is one has
\begin{equation}
\label{E1}
n(f)_a=n(f)_b
\end{equation}
for any integers $a,b$ satisfying $a+b=T$, see \cite{Se, SW}. It was shown in \cite[Corollary 4.3]{DPop} that the graded $S$-module  $N(f)$ satisfies a Lefschetz type property with respect to multiplication by generic linear forms. This implies in particular the inequalities
\begin{equation}
\label{in} 
0 \leq n(f)_0 \leq n(f)_1 \leq ...\leq n(f)_{[T/2]} \geq n(f)_{[T/2]+1} \geq ...\geq n(f)_T \geq 0.
\end{equation}
Moreover, for a degree $d$ 3-syzygy curve $C$ with exponents $(d_1,d_2,d_3)$, we have the following formula for the initial degree of the graded module $N(f)$, see \cite[Theorem 3.9]{minTjurina}.
\begin{equation}
\label{ID}
\sigma (C)= \min \{k \ : \ N(f)_k \ne 0 \}=3(d-1)-(d_1+d_2+d_3).
\end{equation}
For a reduced curve $C$, we denote by $\tau(C)$ its total Tjurina number, that is the sum of the Tjurina numbers of all the singularities of $C$. The following result shows that the invariants $\nu(C)$ and $\tau(C)$ are closely related.
\begin{thm}
\label{thmN}
Let $C:f=0$ be a reduced plane curve of degree $d$ in $\PP^2$ and let $r={\rm mdr}(f)$.
Then the following hold.
\begin{enumerate}
\item If $r < d/2$, then
$$\nu(C)=(d-1)^2-r(d-1-r)-\tau(C).$$

\item If $r \geq (d-2)/2$, then 
$$\nu(C)= \ceil*{ \frac{3}{4}(d-1)^2 } -\tau (C).$$

\end{enumerate}
\end{thm}
Here, for any real number $u$,  $\lceil     u    \rceil $ denotes the round up of $u$, namely the smallest integer $U$ such that $U \geq u$.
Written down explicitly, this means that for $d=2m$ even and $r \geq m-1$, one has
$\nu(C)=3m^2-3m+1-\tau(C)$, while for $d=2m+1$  odd and  $r \geq m$, one has
$\nu(C)=3m^2-\tau(C).$ For $(d-2)/2 \leq r <d/2$, both formulas in (1) and (2) apply, and they give the same result for $\nu(C)$.

\medskip

We continue this section by recalling the following result due to  du Plessis and Wall, see \cite[Theorem 3.2]{duPCTC} as well as \cite{E} for an alternative approach.
\begin{thm}
\label{thmCTC}
For positive integers $d$ and $r$, define two new integers by 
$$\tau(d,r)_{min}=(d-1)(d-r-1)  \text{ and } 
\tau(d,r)_{max}= (d-1)^2-r(d-r-1).$$ 
Then, if $C:f=0$ is a reduced curve of degree $d$ in $\PP^2$ and  $r={\rm mdr}(f)$,  one has
$$\tau(d,r)_{min} \leq \tau(C) \leq \tau(d,r)_{max}.$$
Moreover, for $r={\rm mdr}(f) \geq d/2$, the stronger inequality
$$\tau(C) \leq \tau(d,r)_{max} - \binom{2r+2-d}{2}$$
holds.

\end{thm}

\begin{rk}
\label{rkCTC} Let $C:f=0$ be a reduced curve of degree $d$ in $\PP^2$ and  $r={\rm mdr}(f)$. Note that the function $\tau(d,r)_{max}$, regarded as a function of $r$, occurs also in Theorem \ref{thmN} (1), which can be restated as
\begin{equation}
\label{eqN}
\tau(C)+\nu(C)=\tau(d,r)_{max},
\end{equation}
for $r={\rm mdr(f)} <d/2$.
The inequality
$\tau(C) \leq \tau(d,r)_{max}$ in Theorem \ref{thmCTC} is made more precise, when $r <d/2$, by the result in Theorem \ref{thmN} (1).

\end{rk}

At the end of the proof of Theorem \ref{thmCTC}, in \cite{duPCTC}, the authors state the following very interesting consequence (of the proof, not of the statement) of Theorem \ref{thmCTC}.
\begin{cor}
\label{corCTC} Let $C:f=0$ be a reduced curve of degree $d$ in $\PP^2$ and  $r={\rm mdr}(f)$. One has
$$ \tau(C) =\tau(d,r)_{max}$$
if and only if $C:f=0$ is a free curve, and then $r <d/2$.
\end{cor}
Since a plane curve $C$ is free if and only if $\nu(C)=0$, this characterization of free curves follows also from Theorem \ref{thmN}, as explained in Remark \ref{rkCTC}.

In the paper \cite{Dmax}, we have given an alternative proof of Corollary \ref{corCTC} and have
shown that  a plane curve $C$ is nearly free, which can be defined by the property $\nu(C)=1$, if and only if a similar property holds. Namely, one has the following result, an obvious consequence of Theorem \ref{thmN} and Theorem \ref{thmCTC}. 

\begin{cor}
\label{corCTC2} Let $C:f=0$ be a reduced curve of degree $d$ in $\PP^2$ and  $r={\rm mdr}(f)$.
One has
$$ \tau(C) =\tau(d,r)_{max}-1$$
if and only if $C:f=0$ is a nearly free curve, and then $r  \leq d/2$.
\end{cor}

\section{Terao's Conjecture and some generalizations } 

We say that Terao's Conjecture  holds  for a free line arrangement $\A$ if any other line arrangement $\B$  having an  isomorphic intersection lattice  $L(\B)=L(\A)$,  is also free, see \cite{DHA,OT,Yo14}. 
Terao's Conjecture is known to hold in many cases, for instance when the number of lines in $\A$ is at most 14, see \cite{BK,DIM13}. The open question is the following.
\begin{conj}
\label{conj0}
Terao's Conjecture  holds  for any free line arrangement $\A$.
\end{conj}

\begin{rk}
\label{rkconj0}
Since the total Tjurina number $\tau(\A)$ is determined by the intersection lattice $L(\A)$, see formula \eqref{mu} below, a possible approach to proving Terao's Conjecture may be  to check that $\A:f=0$ and $\B:g=0$ satisfy ${\rm mdr}(f)={\rm mdr}(g)$ and then apply Corollary  \ref{corCTC}. However,
when $\A:f=0$ is a non-free line arrangement, examples due to G. Ziegler in \cite{Z} show that 
the invariant ${\rm mdr}(f)$ is not combinatorially determined.  Consider the following two arrangements 
$$\A: f=xy(x-y-z)(x-y+z)(2x+y-2z)(x+3y-3z)(3x+2y+3z)$$
$$(x+5y+5z)(7x-4y-z)=0,$$
and respectively by
$$\A': f'=xy(4x-5y-5z)(x-y+z)(16x+13y-20z)(x+3y-3z)(3x+2y+3z)$$$$(x+5y+5z)(7x-4y-z)=0$$
see  \cite[Remark 8.5]{DHA}, but beware a misprint in the equation for $\A'$ given there.
This pair of arrangements  satisfy ${\rm mdr}(f)=5$ and ${\rm mdr}(f')=6$, though $\A$ and $\A '$ have the same combinatorics. More precisely, the exponents of $\A$ are $(5,6,6,6)$, while  the exponents of $\A'$ are $(6,6,6,6,6,6)$.

Both arrangements $\A$ and $\A'$ have 6 triple points and 18 double points. For the arrangement $\A$, the 6 triple points are situated on a conic, while for the arrangement $\A'$ this is not the case. An insight into the geometry of the arrangement $\A$ is provided in the paper \cite{DZ}.
\end{rk}

\begin{rk}
\label{rkZ}
There is an interesting  question to find new pairs of line arrangements,
say $\B:g=0$ and $\B':g'=0$ such that their intersection lattices verify $L(\B)=L(\B') \ne L(\A)$ and
${\rm mdr}(g)\ne {\rm mdr}(g').$ 
One way to do this is to add lines to the line arrangements $\A$ and $\A'$ in Remark \ref{rkconj0}.  Add a generic line $L_1$  to the arrangement $\A$
and get in this way a new arrangement $\A_1:f_1=0$. Then add a generic line $L_1'$  to the arrangement $\A'$
and get in this way a new arrangement $\A_1':f_1'=0$. Using \cite[Theorem 3.3]{ADS} or \cite[Corollary 6.4]{DIS}, we see that ${\rm mdr}(f_1)=6$ and ${\rm mdr}(f_1')=7$, though clearly $\A_1$ and $\A _1'$ have the same combinatorics. Then choose a double point $p \in \A_1$ and a double point $p' \in \A_1'$, which correspond to each other under the isomorphism of intersection lattices $L(\A_1)=L(\A_1')$.
Then add generic lines $L_2$ passing through $p$ and $L_2'$ passing through $p'$ to $\A_1$ and
resp.  to $\A_1'$ to get new arrangements $\A_2:f_2=0$ and $\A_2':f_2'=0$.
Using \cite[Theorem 3.3]{ADS} or \cite[Corollary 6.4]{DIS}, we see that ${\rm mdr}(f_2)=7$ and ${\rm mdr}(f_2')=8$, though clearly $\A_2$ and $\A _2'$ have the same combinatorics. Finally, choose a double point $q \in \A_2$ and a double point $q' \in \A_2'$, which correspond to each other under the isomorphism of intersection lattices $L(\A_2)=L(\A_2')$. Add a generic line $L_3$ through the point $q$ to the arrangement $\A_2$
and get in this way a new arrangement $\A_3:f_3=0$. Then add a generic line $L_3'$ through the point $q'$ to the arrangement $\A_2'$
and get in this way a new arrangement $\A_3':f_3'=0$. Using \cite[Theorem 3.3]{ADS} or \cite[Corollary 6.4]{DIS}, we see that ${\rm mdr}(f_3)=8$ and ${\rm mdr}(f_3')=9$, though clearly $\A_1$ and $\A _1'$ have the same combinatorics. In particular, they have both 8 triple points and 42 nodes. Note that the first two steps of the above construction can be interchanged and that by varying the choices of the points $p$ and $q$ several isomorphism classes of intersection lattices
$L(\A_2)$ and $L(\A_3)$ can be obtained.

A different approach for the construction of pairs  $\B:g=0$ and $\B':g'=0$ such that their intersection lattices verify $L(\B)=L(\B')$ and
${\rm mdr}(g)\ne {\rm mdr}(g')$ can be found in \cite{DST}.
\end{rk}

\begin{rk}
\label{rkconj01}
Note that if the line arrangement $\A:f=0$ is free and the line arrangement $\B:g=0$ satisfies $L(\A)=L(\B)$, this clearly implies
$\deg(f)=\deg(g)$ and also ${\rm mdr}(g)\leq {\rm mdr}(f)$ in view of Theorem \ref{thmCTC} and Corollary \ref{corCTC}. It follows that Terao's Conjecture is implied by the following.
\end{rk}
\begin{conj}
\label{conj001}
Let $\A:f=0$ be a  line arrangement in $\PP^2$ with $d=\deg f$ and ${\rm mdr}(f)<d/2$. Then the  integer ${\rm mdr}(f)$ is combinatorially determined.
\end{conj}

\begin{rk}
\label{rkconj011}
We say that two arrangements $\A:f=0$ and $\B:g=0$ have the same weak combinatorial type if for any integer $k \geq 2$, the number of points of multiplicity $k$ is the same in both $\A$ and $\B$. There are examples of pairs of arrangements $\A$ and $\B$, having the same weak combinatorial type, and such that $\A$ is free and $\B$ is not free, see \cite{MaVa2}. Since the  weak combinatorial type determines the total Tjurina number by the formula \eqref{mu}, it follows that in these examples one has ${\rm mdr}(g)< {\rm mdr}(f)<d/2$.
\end{rk}

Theorem \ref{thmN} suggests that the following stronger version of  H. Terao's Conjecture \ref{conj0} might be true, which is clearly equivalent to Conjecture \ref{conj001}.

\begin{conj}
\label{conj1}
Let $\A:f=0$ be a line arrangement in $\PP^2$. Then the invariant $\nu(\A)$ is combinatorially determined.
\end{conj}

For a reduced plane curve one may state the following.
\begin{conj}
\label{conj1G}
Let $C:f=0$ be a reduced plane curve $\PP^2$. Then the invariant $\nu(C)$ is determined by the degree of $C$ and the list of the analytic types of the isolated singularities of $C$.
\end{conj}
Note that both these conjectures hold when $r={\rm mdr}(f) \geq (d-2)/2$ in view of Theorem \ref{thmN}. In particular, the two arrangements $\A$ and $\A'$ in Remark \ref{rkconj0} satisfy $\nu(\A)=\nu(\A')=6$.

\medskip

Consider the sheafification 
$E_C:= \widetilde{D_0(f)} $
of the graded $S$-module $D_0(f)$, which is a rank two vector bundle on $\PP^2$, see \cite{Se} for details. Moreover, recall that
\begin{equation} \label{equa1} 
E_C=T\langle C \rangle (-1),
\end{equation}
where $T\langle C \rangle $ is the sheaf of logarithmic vector fields along $C$ as considered for instance in \cite{DSer,MaVa,Se}. 
One has,  for any integer $k$, 
\begin{equation}
\label{e7}
H^0(\PP^2, E_C(k))=D_0(f)_k \text{ and }  H^1(\PP^2, E_C(k))=N(f)_{k+d-1},  
\end{equation}
where $d=\deg(f)$, for which we refer 
to \cite[Proposition 2.1]{Se}. Note that $C:f=0$ is a free curve with exponents $(d_1,d_2)$ if and only if the vector bundle $E_C$ splits as a direct sum
\begin{equation} \label{equa1a} 
E_C=\OO_{\PP^2}(-d_1) \oplus \OO_{\PP^2}(-d_2).
\end{equation}
In general, let  $(d_1^{L_0},d_2^{L_0})$, with  $d_1^{L_0} \leq d_2^{L_0}$, denote the splitting type of the vector bundle $E_C$ along a generic line $L_0$ in $\PP^2$.
Then Proposition 3.2 in \cite{AD} implies that 
$d_1^{L_0}=r$ for $r <(d-2)/2$ and $d_1^{L_0}=\lfloor (d-1)/2 \rfloor $ for $r  \geq (d-2)/2$, where $r={\rm mdr}(f)$.  Then \cite[Theorem 1.1]{AD} says that 
$$(d-1)^2 -d_1^{L_0}d_2^{L_0}=\tau(C)+ \nu(C).$$
Since $d_1^{L_0}+d_2^{L_0}=d-1$ by \cite[Proposition 3.1]{AD},  Conjecture \ref{conj1} may be restated as follows.

\begin{conj}
\label{conj2}
Let $C:f=0$ be a line arrangement in $\PP^2$. Then the generic splitting type  $(d_1^{L_0},d_2^{L_0})$ of the vector bundle $E_C=T\langle C\rangle(-1)$ is combinatorially determined.
\end{conj}
We note that Conjecture \ref{conj2} is just the question asked in \cite[Question 7.12]{Cook+}.

\section{Some related results} 
Here we discuss an alternative view on the $S$-module
$$D_0(f)=AR(f)={\rm Syz}(f).$$
Let $\Omega^j$ by the graded $S$-module of global, polynomial differential $j$-forms in $\C^3$, for $j=0,1,2,3$. Note that $\Omega^0=S$ and $\Omega^3=S \dd x \wedge \dd y \wedge \dd z=S(-3)$.
To a triple $\rho=(a,b,c) \in S^3$ we associate the 2-form
$$\omega(\rho)=a\dd y \wedge \dd z -b \dd x \wedge  \dd z +c \dd x \wedge  \dd y.$$
Then $\rho \in D_0(f)$ if and only if $ \dd f \wedge \omega(\rho)=0$,
where $\dd f=f_x \dd x +f_y \dd y +f_z \dd z$ is the differential of $f$. In other words, we have
\begin{equation} \label{eqAR} 
AR(f)(-2)= \ker \{ \dd f : \Omega ^2 \to \Omega ^3\}.
\end{equation}
We define the submodule of Koszul-type relations $KR(f)$ to be
\begin{equation} \label{eqKR} 
KR(f)(-2)= \im \{ \dd f : \Omega ^1 \to \Omega ^2\}
\end{equation}
and note that $KR(f)$ is generated by 3 obvious relations of degree $d-1$, namely
$$(f_y,-f_x,0), \  (f_z,0,-f_x)  \text{  and  } (0,f_z,-f_y).$$
Finally, consider the quotient of essential relations
\begin{equation} \label{eqKR2} 
ER(f)= AR(f)/KR(f),
\end{equation}
or, in cohomological terms, $ER(f)(-2)=H^2(\Omega^*, \dd f \wedge)$.
To state the following key result,  we recall some more notation. 
Let $J=J_f$ be the Jacobian ideal of $f$ and $I=I_f$ be its saturation with respect to the maximal ideal $(x,y,z)$. Then the singular subscheme $\Sigma_f$ of the reduced curve $C:f=0$ is the 0-dimensional scheme defined by the ideal $I$ and we consider the following sequence of defects
\begin{equation} 
\label{eqDEF}
\defect _k\Sigma_f=\tau(C)-\dim \frac{S_k}{I_k}.
\end{equation}
With this notation, one has the following result, see \cite[Theorem 1]{Bull13}.

\begin{thm}
\label{linsys}
Let $C:f=0$ be a degree $d$ reduced curve in $\PP^2$. If $\Sigma_f$ denotes its singular locus subscheme, then
$$\dim ER(f)_{2d-5-k}= \defect _k\Sigma_f $$
for $0\leq k \leq 2d-5$  and $\dim ER(f)_j=\tau(C)$ for $j\geq 2d-4$.
\end{thm}
It follows that Terao's Conjecture \ref{conj0} is equivalent to the following, recall Conjecture \ref{conj001}.
\begin{conj}
\label{conj00}
Let $\A:f=0$ be a  line arrangement in $\PP^2$ with $d=\deg f$ and ${\rm mdr}(f)<d/2$. Then the largest integer $k$ such that   $\defect _k\Sigma_f \ne 0$ is combinatorially determined.
\end{conj}

\begin{cor}
\label{corlinsys}
Let $C:f=0$ be a degree $d$ free curve in $\PP^2$, with exponents $(d_1,d_2)$. Then
$$\defect _k\Sigma_f=\tau(C)-\dim M(f)_k=0$$
for $k >k_1= 2d-5-d_1$  and
$$\defect _{k_1}\Sigma_f=\tau(C)-\dim M(f)_{k_1}>0.$$
\end{cor}

Note that we have
\begin{equation} 
\label{eqDEF1}
\defect _k\Sigma_f=\tau(C)-\dim \frac{S_k}{J_k}+ n(f)_k.
\end{equation}
We introduce the following.
\begin{definition}
\label{defCTST} For a homogeneous reduced polynomial $f \in S_d$ one defines
\begin{enumerate}
\item [(i)] the {\it coincidence threshold}\index{coincidence threshold ${\rm ct}(f)$}
$${\rm ct}(f)=\max \{q:\dim M(f)_k=\dim M(f_s)_k \text{ for all } k \leq q\},$$
with $f_s$  a homogeneous polynomial in $S$ of the same degree $d$ as $f$ and such that $C_s:f_s=0$ is a smooth curve in $\PP^2$.
\item [(ii)] the {\it stability threshold}
$${\rm st}(f)=\min \{q~~:~~\dim M(f)_k=\tau(C) \text{ for all } k \geq q\}.$$
\end{enumerate}

\end{definition}
It is clear that one has
\begin{equation}
\label{REL}
{\rm ct}(f) \geq {\rm mdr}(f)+d-2,
\end{equation}
with equality for ${\rm mdr}(f) <d-1.$ The invariant ${\rm st}(f)$ is more mysterious, the equality in \eqref{eqDEF1} implies that
$${\rm st}(f)=\min \{q~~:~~\dim N(f)_k=\defect_k \Sigma_f \text{ for all } k \geq q\}.$$
These new invariants ${\rm ct}(f)$ and ${\rm st}(f)$ enter into the following result, see \cite[Corollary 1.7]{Dmax} and recall that $T=3(d-2)$.
\begin{thm}
\label{ct+st}
Let $C:f=0$ be a degree $d$ reduced curve in $\PP^2$. Then $C$ is a free (resp. nearly free) curve if and only if 
$${\rm ct}(f)+ {\rm st}(f)=T \text{ (resp.  }   \  {\rm ct}(f)+ {\rm st}(f)=T+2). $$
In the remaining cases one has ${\rm ct}(f)+ {\rm st}(f)\geq T+3.$
\end{thm}
It follows that a new stronger form of Terao's Conjecture might be the following.

\begin{question}
\label{qT4}
Let $\A:f=0$ be a line arrangement in $\PP^2$. Is the invariant 
$${\rm ct}(f)+ {\rm st}(f)$$
  combinatorially determined ?
\end{question}
 For the two degree 9 line arrangements considered by Ziegler and recalled above in Remark \ref{rkconj0}, we have ${\rm ct}(f)=12$,  ${\rm st}(f)=14$, ${\rm ct}(f')=13$ and   ${\rm st}(f')=13$, hence 
$${\rm ct}(f)+ {\rm st}(f)={\rm ct}(f')+ {\rm st}(f')= T+5=26  \text{  and  }  \nu(\A)=6.$$
However, for the pair of line arrangements $\A_1$ and $\A_1'$, constructed as in Remark \ref{rkZ}, one has
 ${\rm ct}(f_1)=14$,  ${\rm st}(f_1)=16$ and  ${\rm ct}(f_1')=15$,  ${\rm st}(f_1')=16$. This shows that Question \ref{qT4} has a negative answer.

 The relation between this invariant ${\rm ct}(f)+ {\rm st}(f)$ and the defect of freeness $\nu(C)$ is also an {\it open problem}.
For a plus-one generated curve $C:f=0$ we have
$${\rm ct}(f)+ {\rm st}(f)=T+ \nu(C)+1,$$
 see \cite[Proposition 3.7]{minTjurina}.
 However, the formulas in \eqref{nuVS} below show that such a linear dependence cannot hold in general, even for 3-syzygy curves.

\section{Freeness properties of rational cuspidal curves } 

A plane rational cuspidal curve is a rational curve $C:f=0$ in the complex projective plane $\PP^2$, having only unibranch singularities. The study of these curves has a long and fascinating history, some long standing conjectures, as the Coolidge-Nagata conjecture being
proved only recently, see \cite{KP0}.
The classification of such curves is not easy, there are a wealth of examples even when additional strong restrictions are imposed, see \cite{FLMN, FZ, Moe,PP}.

We have remarked in \cite{MRL} that many plane rational cuspidal curves are free. The remaining examples of plane rational cuspidal curves in the available classification lists turned out to satisfy a weaker homological property, which was chosen as the definition of a nearly free curve, see \cite{DStRIMS}. Subsequently, a number of authors have establish interesting properties of this class of curves, see \cite{B+, MaVa}. In view of the above remark, we have conjectured in \cite[Conjecture 1.1]{DStRIMS}  the following surprising fact.

\begin{conj}
\label{conjRCC}
Any plane rational cuspidal curve $C$ is either free or nearly free. 
\end{conj}
This conjecture was proved in \cite[Theorem 3.1]{DStRIMS} for curves $C$ whose degree $d$ is even, as well as for some cases when $d$ is odd, e.g. when $d=p^k$, for a prime number $p>2$. It turns out that a rational cuspidal curve $C:f=0$ with ${\rm mdr}(f)=1$ is nearly free. Indeed,  this follows from \cite[Proposition 4.1]{Drcc}. To see this, note that the implication (1) $\Rightarrow $ (2) there  holds for any $d \geq 2$.

Assume from now on that $d$ is odd, and let 
\begin{equation}
\label{pfactor}
d=p_1^{k_1} \cdot p_2^{k_2} \cdots p_m^{k_m}
\end{equation}
be the prime decomposition of $d$. We assume also that $m\geq 2$, the case $m=1$ of our conjecture being settled in \cite[Corollary 3.2]{DStRIMS}. By changing the order of the primes $p_j$'s if necessary, we can and do assume that  $p_1^{k_1}>p_j^{k_j}$,
for any $2 \leq j \leq m$. Set $e_1=d/p_1^{k_1}.$
With these assumptions and notations, we have  the following results, see \cite{Mosk}.

\begin{thm}
\label{thmA5}
Let $C:f=0$ be a rational cuspidal curve of degree $d=2d'+1$ an odd number. Then ${\rm mdr}(f) \leq d'$ and if equality holds, then $C$ is either free or nearly free.

\end{thm} 

\begin{thm}
\label{thmB5}
Let $C:f=0$ be a rational cuspidal curve of degree $d=2d'+1$,  an odd number as in \eqref{pfactor}. Then, if
$${\rm mdr} (f) \leq r_0:=\frac{d-e_1}{2},$$ then $C$ is either free or nearly free. In particular, the following hold.

\begin{itemize}
		\item[i)]  If $d=3p^k$, with $p$ a prime number, then $C$ is either free or nearly free.
		\item[ii)]  $d=5p^k$, with $p$ a prime number, $p^k>3$, then $C$ is either free or nearly free, unless ${\rm mdr}(f)=d'-1$.
		
	\end{itemize}

\end{thm}

\begin{rk}
\label{rkA5}
Note that, for $d \ne 15$, we have $e_1 \leq d/7$ and hence
$$r_0=\frac{d-e_1}{2} \geq \ceil*{\frac{d(1-\frac{1}{7})}{2}}=\ceil*{\frac{3d}{7}}.$$
Therefore, the only cases not covered by our results correspond  to curves of odd degree $d$, such that $r={\rm mdr}(f)$ satisfies
$$\ceil*{\frac{3d}{7}}+1 \leq r_0+1 \leq r \leq d'-1=\frac{d-3}{2}.$$
\end{rk}

\begin{cor}
\label{corA5}
A rational cuspidal curve $C:f=0$ of degree $d$ with $r={\rm mdr}(f)$ is either free or nearly free, if 
one of the following conditions holds.

\begin{enumerate}

\item $r \leq 15$, or

\item $d \leq 90$, unless we are in one of the following situations.

\begin{itemize}
		\item[i)]  $d=35$ and $r=16$;
		\item[ii)]  $d=45$ and $r=21$;
		\item[iii)]  $d=55$ and $r=26$;
		\item[iv)]  $d=63$ and $r \in \{29,30\}$;
		\item[v)]  $d=65$ and $r=31$.
		\item[vi)]  $d=77$ and $r \in \{36,37\}$.
		\item[vii)]  $d=85$ and $r=41$.
	\end{itemize}

\end{enumerate}

\end{cor}
In the excluded situations, our results do not allow us to conclude, and hence Conjecture \ref{conjRCC} is still open.

Here is a recent additional result in this area, see \cite[Theorem 1.5]{POG}. Consider an affine plane curve $X:g(x,y)=0$ given by a reduced polynomial $g \in R=\C[x,y]$ of degree $d$. Then the projective closure
$\overline X$ of $X$ is the curve in $\PP^2$ defined by the polynomial
$$f(x,y,z)=z^dg(\frac{x}{z},\frac{y}{z}).$$
Recall that  a contractible, irreducible affine plane curve $X$ is given, in a {\it suitable global coordinate system} $(u,v)$ on $\C^2$ by the equation $u^p-v^q=0$ for some relatively prime integers $p \geq 1$ and $q \geq 1$, see \cite{GM,LZ}. In particular, $X$ has at most a unique singular point $a$, which is a cusp of type $(p,q)$, namely the singularity $(X,a)$ is given in local analytic coordinates $(u',v')$ at $a$ by the equation $u^{'p}-v^{'q}=0$. When $X$ is smooth, then $X$ is isomorphic to $\C$ and $g$ is a component of an automorphism of $\C^2$, see
\cite{Ab,Su}. We say that in this case $X$ has a cusp of type $(1,1)$.

\begin{thm}
\label{thmC5}
With the above notation, assume that $X$ is irreducible and contractible and has a cusp of type $(p,q)$ such that either $p$ or $q$ is relatively prime to $d+1$, where $d= \deg X= \deg \overline X$.
Then the  projective closure $\overline X$ of the affine plane curve $X$ is a rational cuspidal curve which is either free or plus-one generated.
\end{thm}
All the proofs of the results in this section, in spite of their purely algebraic statements, use a deep result due to U. Walther involving mixed Hodge theory and D-modules, see  \cite[Theorem 1.6]{Wa}.

\section{Curves with minimal and maximal total Tjurina numbers}

Consider a reduced curve $C:f=0$ of degree $d$ and set $r={\rm mdr}(f)$.
Then Theorem \ref{thmCTC} gives upper and lower bounds for the total Tjurina number $\tau(C)$ as functions of $d$ and $r$.
Moreover, we know that when $r<d/2$ the upper bound is obtained exactly for a free curve, recall Corollary \ref{corCTC}. Moreover, it was shown in \cite{expo} that for any $r<d/2$ there is a free curve with exponents $d_1=r$ and $d_2=d-1-r$. The situation of the equality
\begin{equation} 
\label{eqMAX}
\tau(C) = \tau(d,r)_{max} - \binom{2r+2-d}{2}
\end{equation}
for 
$r \geq d/2$ seems to be more subtle.  A reduced curve $C:f=0$
of degree $d$, such that $r={\rm mdr} (f)$ and the equality \eqref{eqMAX} holds for $C$ is called a {\it maximal Tjurina curve of type} $(d,r)$.
In the paper \cite{maxTjurina} we have put forth the following.
\begin{conj}
\label{conjMAX}
For any integer $d\geq 3$ and for any integer $r$ such that $d/2\leq  r \leq d-1 $, there
are maximal Tjurina curves of type $(d, r)$. Moreover, for $d/2\leq  r \leq d -2$, there are maximal
Tjurina line arrangements of type $(d, r)$.
\end{conj}
This conjecture was shown to hold for many pairs $(d,r)$, and in particular for all pairs as above with $d \leq 11$, see \cite{ADS,maxTjurina}, but the general case is still open. In many cases, these constructed maximal Tjurina curves are either line arrangements,
or irreducible curves with interesting geometrical properties, e.g. for the pair $(d,r)=(d,d-1)$ we get as maximal Tjurina curves the maximal nodal curves, namely irreducible nodal curves of degree $d$
having 
$$g=(d-1)(d-2)/2$$
 nodes, see \cite{maxTjurina, Oka}. 
 
 One can asks for upper bounds for the total Tjurina number of a hypersurface $V$ in $\PP^n$ having only isolated singularties, in the case $n \geq 3$. In this setting, some results and some {\it open problems} are stated in \cite{DAndr}.
 
  \medskip
 
 Let's have a look now at the curves $C$ for which the total Tjurina number $\tau(C)$ has the minimal possible value given in Theorem \ref{thmCTC}. The Thom-Sebastiani plane curves, that is the plane curves 
 $$C: f(x,y,z)=g(x,y)+z^d=0$$
 where $g(x,y)$ is a homogeneous polynomial in $x$ and $y$ of degree $d$ show that this lower bound is effective. In other words they give examples of curves for any $d$ and $r={\rm mdr}(f)$ such that
 $$\tau(C)=\tau(d,r)_{min}=(d-1)(d-r-1),$$
 see \cite[Example 4.5]{minTjurina}. More precisely, if $m$ denotes the number of distinct linear factors of $g(x,y)$, then $r=m-1$ and
 $\tau(C)=(d-1)(d-m)$. With this notation, one also has
\begin{equation}
\label{nuVS}
{\rm ct}(f)+ {\rm st}(f)=(d+m-3)+(2d+m-5)=T+2m-2,
\end{equation} 
while $\nu(C)=m(m-2)+1$, see  \cite[Example 3.8]{minTjurina}.
This shows that the answer to the open problem mentioned after Question \ref{qT4} is not easy.
 \medskip
 
 It is a remarkable fact in our opinion that this lower bound  $\tau(d,r)_{min}$ is not optimal if we restrict our attention to line arrangements, as shown in \cite{minTline}. 
 When $C$ is a line arrangement,  its global Tjurina number, coincides with its global Milnor number $\mu(C)$, and is given by
\begin{equation}
\label{mu}
\tau(C)=\sum_p(n(p)-1)^2,
\end{equation}
the sum being over all multiple points $p$ of $C$, and $n(p)$ denoting the multiplicity of $C$ at $p$.
Let $m(C)$ be the maximal multiplicity of a point in $C$,
and $n(C)$ the maximal multiplicity of a point in $C \setminus \{p\}$, where $p$ is any point in $C$ of multiplicity $m(C)$. Note that 
$$ 1 \leq n(C) \leq m(C) \leq d.$$
Moreover $m(C)=d$ if and only if ${\rm mdr}(f)=0$, and $m(C)=d-1$ if and only if ${\rm mdr}(f)=1$,
see \cite[Proposition 4.7]{DIM}. 
With this notation we have the following result, see \cite{minTline}. 

\begin{thm}
\label{thmA6}
Let $C:f=0$ be an arrangement of $d\geq 4$ lines in $\PP^2$ which is not free. If we set $r={\rm mdr}(f)\geq 2 $ and $\tau(d,r)_{min}=(d-1)(d-r-1)$, then  the following hold.

\begin{itemize}
		\item[i)]  With the above notation, one has
		$$\tau(C) \geq \tau'(d,r)_{min}:= \tau(d,r)_{min}+\binom{r}{2}+\binom{n(C)}{2}+1.$$
\item[ii)] If $r \ne d-m(C)$, then the possibly stronger inequality 	
$$\tau(C) \geq \tau''(d,r)_{min}:=\tau(d,r)_{min}+\binom{r}{2} +\binom{m(C)}{2}+1$$ holds.	

\end{itemize}

\end{thm}
The line arrangements such that $r={\rm mdr}(f) \in \{0,1,2\}$ are classified, see \cite{T2} or \cite[Theorem 4.11]{DIM} for the case $r=2$, which is the only difficult case. In addition, the case
$2=n(C)\leq m(C) \leq d-2$,  corresponds to the intersection lattice $L(C)$ being the lattice $L(d,m(C))$ discussed in \cite[Proposition 4.7]{DIM}, i.e. the intersection lattice of an arrangement obtained from $m(C)$ concurrent lines by adding $d-m(C)$ lines in general position.
For the remaining line arrangements we have the following.

\begin{cor}
\label{corA6} Let $C:f=0$ be an arrangement of $d$ lines in $\PP^2$ which is not free and such that $r={\rm mdr}(f)\geq 3 $ and $n(C) \geq 3$. Then 
$$\tau(C) \geq \tau^N(d,r)_{min}:=\tau(d,r)_{min}+\binom{r}{2} +4 \geq \tau(d,r)_{min}+7.$$
\end{cor}
However, it is an {\it open question to find the
 best lower bound for the global Tjurina number of line arrangements}.
 
 \section{Free curves with ADE singularities}
 As we have already remarked above,  it was shown in \cite{expo} that for any $r<d/2$ there is a free  curve $C$ with exponents $d_1=r$ and $d_2=d-1-r$. We can ask about the existence of such a curve when we impose the extra condition that $C$ has only simple singularities $A,D,E$. This question arose naturally in the paper
\cite{DPok}, where we showed that the maximizing curves introduced by U. Persson \cite{Persson} were in fact a class of free curves with $A,D,E$ singularities.
Here are the necessary definition. 
\begin{definition}
\label{defmax7}
A reduced curve $C$ in $\PP^2$ of even degree $d=2m\geq 4$ is a {maximizing curve} if $C$ has only simple singularities $A,D,E$ and
$$\tau(C) = 3m(m-1)+1.$$
\end{definition}
For such a maximizing curve $C$, the associated smooth surface $\widetilde{X}$ obtained as the minimal resolution of the double cover $X$ of $\PP^2$ ramified along $C$ might have the maximal possible Picard number. From this perspective, it is very interesting to classify, or at least to produce  many examples of such maximizing curves. The main result of the paper \cite{DPok} can be formulated as follows.
\begin{thm}
\label{maxi7}
Let $C$ be a plane curve of degree $d=2m \geq 4$ having only ${A,D,E}$ singularities. Then $C$ is a maximizing curve if and only if $C$ is a free curve with the exponents $(m-1,m)$.
\end{thm}
Moreover, it was shown in \cite{Persson} that maximazing curves exist for any even degree $d=2m \geq 4$. For instance the family of curves
in \cite[Lemma 7.8]{Persson} 
$$\mathcal C_{2(m+1)}: xy[(x^m+y^m+z^m)^2-4(x^my^m+y^mz^m+z^mx^m)]=0$$
is maximizing, and hence we get a curve of even degree $d=2(m+1)$ which is free with exponents $(d_1,d_2)=(m,m+1)$. Motivated by Definition \ref{defmax7} we introduced in \cite{DPok} the following definitions for the odd degree curves.
\begin{definition}
\label{defmaxODD}
A reduced  curve $C$ in $\PP^2$ with only $A,D,E$ singularities and of odd degree $d=2m+1\geq 5$ is a {maximizing curve} if
$$\tau(C) = 3m^2+1.$$
\end{definition}
This definition was chosen in view of the following result.
\begin{prop}
\label{propODD1}
Let $C: f=0$ be a reduced curve of degree $n=2m+1\geq 5$ with at most ${ A, D, E}$ singularities.  Then one of the following three situations occurs, where $r={\rm mdr}(f)$.
\begin{enumerate}
\item[a)] $r =m-1$ and  $\tau(C) \leq \tau(n,r)_{max}= 3m^2+1$.

\item[b)] $r =m$ and  $\tau(C) \leq \tau(n,r)_{max}= 3m^2$.

\item[c)] $r >m$ and  $\tau(C) < 3m^2-1$.

\end{enumerate}
In both cases $a)$ and $b)$, the equality holds if and only if the curve $C$ is free.
\end{prop}
The following question is natural.

\begin{question}
\label{qmaxODD} Do maximizing curves $C_{2m+1}$ exist in any degree $d=2m+1\geq 5$ ?
\end{question}

In \cite{DPok} we have constructed maximizing curves of degree 5 and 7.
A maximizing curve of degree 9 is constructed in \cite[Example 6.4]{DIPS}, as follows.
The curve
$$C: f=xyz\left( (x^2+y^2+z^2)^3-27x^2y^2z^2 \right)=0$$
has  six singularities $E_7$ at the points $(0:1:\pm i)$, $(1:0: \pm i )$ and $ (1:\pm i: 0)$, and seven nodes $A_1$ at the points
$(1:\pm 1: \pm 1)$, $(1:0:0)$, $(0:1:0)$ and $(0:0:1)$.
It follows that
$$\tau(C)=6\cdot 7+7 \cdot 1=49,$$
and this equality implies that $C$ is maximizing by Definition \ref{defmaxODD}. Proposition \ref{propODD1} tells us that $C$ is a free curve with exponents $(3,5)$.

 For odd degree $d=2m+1 \geq 11$ it seems that no such maximazing curve is known at this time.
 
 \section{Freeness properties of supersolvable plane curves} 
 We start with some definitions, see \cite{DIPS}.
 \begin{definition}
\label{defSS}
Given a reduced plane curve $C$, we say that $p \in C$ is a {\it modular point} for $C$ if the central projection 
$$\pi_p : \PP^2 \setminus \{p\} \to \PP^1$$
induces a locally trivial fibration of the complement $M(C)=\PP^2 \setminus C$. We say that a 
reduced plane curve $C$ is {\it supersolvable} if it has at least one modular point.
\end{definition}
The map induced by $\pi_p$ is a locally trivial fibration if and only if for any line $L_p$ passing through $p$ and not an irreducible component of $C$, one has the following properties for the intersection multiplicities of the line $L_p$ and the curve $C$:
$$(C,L_p)_p=\mult_p(C) \text{ and } (C,L_p)_q=1 \text{ for any  } q \in C \cap L_p, \ q \ne p.$$
 When $C$ is a line arrangement, this definition of a modular point coincides with the usual one,  and a line arrangement is supersolvable by definition if it has a modular point. In particular, the existence of a modular point for a line arrangement $C$ implies that $C$ is free,  see for all these well known facts \cite{DHA,OT}.  The following conjecture was stated in \cite[Conjecture 1.10]{DIPS}.
\begin{conj}
\label{conj18}
Any supersolvable plane curve is free.
\end{conj}
One rather general setting where this conjecture holds is the following, see \cite[Theorem 1.11]{DIPS}.

\begin{thm}
\label{thm58}
Let $C_0$ be a reduced plane curve, let $p \in M(C_0)$ be a point in the complement of $C_0$ and let $\A$ be the set of lines $L$ passing through $p$ such that there is a point $q \in L \cap C_0$ with $(C_0,L)_q >1$. Assume that all the singularities of the curve $C$ obtained by adding all the lines in $\A$ to $C_0$ are quasi homogeneous. Then $C$ is supersolvable and free.
In particular, this holds when all the singularities $s \in C_0$ have multiplicity 2, and $p$ is not on any tangent cone $TC_s(C_0)$ for $(C_0,s)$ a singularity of $C_0$ with $\mu(C_0,s) \geq 3$.
\end{thm}
However, many other supersolvable curves are known to be free.
For instance,
the  free curve 
$$C: f=yz(y^m+z^m)(x^my^m+y^mz^m+x^mz^m)=0$$ 
constructed  in \cite[Theorem 1.8]{DIPS} is of a different nature, since in this case 
$$p \in C_0:x^my^m+y^mz^m+x^mz^m=0 $$
 and the curve $C$  has not only quasi homogeneous singularities. Hence it gives new examples where Conjecture \ref{conj18} holds.
In the same paper,  countable examples of conic-line arrangements were constructed, where the Conjecture \ref{conj18} holds even in the presence of non quasi homogeneous singularities, see \cite[Remark 7.4]{DIPS}.
An equation for such a conic-line arrangement is
$$C:f=x(x^m+z^m)(x^{2m}+(xz+y^2)^m)=0.$$
Here $p=(0:1:0)$ is a modular point and $q=(0:0:1) \in C$ is not a quasi homogeneous singularity for $C$.

\end{document}